\input amstex
\documentstyle{amsppt}
\topmatter
\title
Invariants and orthogonal $q$-polynomials associated with  $\Bbb{C}_{q}(osp(1,2))$
\endtitle
\author
Yi Ming Zou
\endauthor
\address
Department of Mathematical Sciences,
University of Wisconsin-Milwaukee,
Milwaukee, WI 53201, USA
\endaddress
\email
ymzou\@uwm.edu
\endemail
\subjclass
17B35, 17B60, 17B70, 22E70
\endsubjclass
\abstract
The spaces of invariants and the zonal spherical functions associated with quantum super 2-shpheres defined by $\Bbb{C}_{q}(osp(1,2))$ are discussed. Connection between the zonal spherical functions and orthogonal $q$-polynomials from the Askey-Wilson scheme is investigated. 
\endabstract
\endtopmatter

\document
\head
1. Introduction
\endhead
The study of the spherical functions on symmetric spaces can be traced back to Cartan in the 1920s.  Since the late 1980s, this theory has been generalized to quantum group and quantum homogeneous spaces began with [Ko, MNU] (see [Di, Le] for an outline of the philosophy of these approaches and the major steps, as well as some recent developments). In this paper, we study homogeneous spaces and spherical functions associated with the quantum supergroup $\Bbb{C}_{q}(osp(1,2))$ introduced in [Z].
\par
Let $G$ be a contragredient Lie superalgebra over the field $\Bbb{C}$ with a symmetrizable Cartan matrix, let $U(G)$ be the universal enveloping algebra of $G$, and let $U_{q}(G)$ be the quantized enveloping algebra of $G$. Both $U(G)$ and $U_{q}(G)$ are $\Bbb{Z}_{2}$-graded Hopf algebras.  Let $\Bbb{C}[G]$ and $\Bbb{C}_{q}[G]$ be the supergroup and the quantum supergroup corresponding to $G$ constructed from the coordinate rings of $U(G)$ and $U_{q}(G)$ respectively.  By using the $\Bbb{Z}_{2}$-graded Hopf algebra structure, we define a left action $L$ and a right action $R$ of $U(G)$ on $\Bbb{C}[G]$ (similarly for the actions of $U_{q}(G)$ on $\Bbb{C}_{q}[G]$) by 
$$
L_{x}(u)=\sum (-1)^{p(x)p(u_{1})}u_{2}(x)u_{1}, \quad
R_{x}(u)=\sum u_{1}(x)u_{2},\tag 1.1
$$
where $x\in U(G)$, $u\in \Bbb{C}[G]$, $p(x)$ is the parity function associated with the $\Bbb{Z}_{2}$-grading, and $\Delta(u) =\sum u_{1}\otimes u_{2}$.  
\par
For $x, y \in U(G)$ and $u, u' \in \Bbb{C}[G]$, we have
$$
\gather
L_{xy} = L_{x}L_{y},\quad R_{xy} = (-1)^{p(x)p(y)}R_{y}R_{x}, \tag 1.2\\
R_{x}L_{y} = (-1)^{p(x)p(y)}L_{y}R_{x}, \tag 1.3\\
L_{x}(uu') = \sum (-1)^{p(u)p(x_{2}})L_{x_{1}}(u)L_{x_{2}}(u'), \tag 1.4 \\
R_{x}(uu') = \sum (-1)^{p(u)p(x_{2}})R_{x_{1}}(u)R_{x_{2}}(u'), \tag 1.5
\endgather
$$
where $x_{1}$ and $x_{2}$ are defined by $\Delta (x) = \sum x_{1}\otimes x_{2}$.
\par
View the elements of $U(G)$ as functions on $\Bbb{C}[G]$, then for $x, y \in U(G)$ and $u \in \Bbb{C}[G]$, we have
$$
(xy)(u) = x(L_{y}(u)), \qquad (xy)(u) = (-1)^{p(x)p(y)}y(R_{x}(u)). \tag 1.6
$$
\par
If $J$ is a coideal (two-sided) of $U(G)$, then the invariant subspace
$$
B_{J} = \{u\in \Bbb{C}[G]: R_{x}(u)=\varepsilon (x)u\},\tag 1.7
$$
where $\varepsilon$ is the counit of $U(G)$, can be considered as the homogeneous space defined by $J$ with the right action of $\Bbb{C}[G]$ given by the comultiplication (see [DK]).  We call the elements of the left $J$-invariant subspace of $B_{J}$
$$
\{u\in B_{J} : L_{x}(u) =\varepsilon (x)u\} \tag 1.8
$$
spherical functions on the homogeneous space $B_{J}$. We remark that by symmetry, one can also first consider the left invariants to obtain homogeneous spaces and then take the right invariants to be the spherical functions.
\par
If the algebra $G$ is a semisimple Lie algebra and the homogeneous space is a compact symmetric space defined by a maximal compact subalgebra $K$ of $G$, then by the Peter-Weyl's theorem, the representative ring of the pair $(G,K)$ is a direct sum of irreducible representations of $G$, and the zonal spherical functions in each of the irreducible components can be expressed by certain special functions (e.g. hypergeometric functions).  It turns out that this theory can be extended to the quantum groups, and the spherical functions associated with the quantum symmetric pairs can be expressed in terms of $q$-hypergeometric polynomials (see [Ko, DN, Le]).  In the super case, since the lack of reducibility for the representations of $G$ in general, one needs to take different approach to the invariant and spherical function theory.  In [ZZ], the first fundamental theorem of invariant theory of general linear supergroup was applied to investigate the spherical functions on certain homogeneous superspaces associated with general linear supergroup, these homogeneous superspaces include the projective superspace $\Bbb{C}[\Bbb{P}^{n-1|m}]$. The results in [ZZ] show that depends on $m$ and $n$, the space of the spherical functions on $\Bbb{C}[\Bbb{P}^{n-1|m}]$ may or may not decomposes into the direct sum of the eigenfunctions of the Laplacian operator (zonal spherical functions). For the contragredient finite dimensional simple Lie superalgebras, $osp(1,2n)$ is the only one (other than the simple Lie algebras) whose finite dimensional representations are completely reducible.  The interests of study $osp(1,2)$ are similar to those of $sl(2)$, since it is a basic building block for contragredient Lie superalgebras. In [Z], a theory parallel to that of the quantum $SL(2)$ was developed for $OSp(1,2)$. There the definition of the quantized enveloping algebra $U_{q}(osp(1,2))$ is standard, while the definition of $\Bbb{C}_{q}(osp(1,2))$ can be explained by the Radford-Majid bosonization. In this paper, we turn our attention to the homogeneous superspaces associated with $\Bbb{C}_{q}(osp(1,2))$ and study the spherical functions on them. 
\par
In section 2, we recall the basic definitions and results associated with the $\Bbb{Z}_{2}$-graded Hopf algebras $\Cal{U}=U_{q}(osp(1,2))$ and $\Cal{A}=\Bbb{C}_{q}(osp(1,2))$.  In section 3, we first classify all $(g,g')$-primitive elements of $\Cal U$, then study the right and the left invariants of these primitive elements.  In section 4, we investigate how the invariants in section 3 are related to the special $q$-polynomials via the Askey-Wilson scheme.  As suggested in [Se], some of the orthogonal polynomials (or $q$-polynomials) obtained from algebras associated with $osp(1,2)$ may be new, but the author has no answer to this question at this time. This question deserves further attention. 
\par
\head
2. Preliminary 
\endhead
\par
We recall the definitions and some basic properties of the quantized eveloping algebra and the quantum supergroup correspond to $osp(1,2)$ from [Z]. We assume through out that $0\ne q\in \Bbb{C}$ and $q$ is not a root of $1$.  Let $\bold{i} = \sqrt{-1}$, and let $t =\bold{i}\sqrt{q}$. The quantized enveloping algebra $\Cal{U}$ of $osp(1,2)$ is the $\Bbb{Z}_{2}$-graded associative algebra over $\Bbb{C}$ generated by $k^{\pm 1}, e, f$ subject to the relations: 
$$\gathered
kk^{-1}=k^{-1}k=1,\\
kek^{-1}=qe, \qquad kfk^{-1}=q^{-1}f,\\
ef + fe = \frac{k-k^{-1}}{q-q^{-1}}.
\endgathered\tag 2.1
$$
The $\Bbb{Z}_{2}$-grading on $\Cal{U}$ is given by $p(k^{\pm 1})=0$ and $p(e)=p(f)=1$.
The $\Bbb{Z}_{2}$-graded Hopf algebra structure on $\Cal{U}$ is defined by
$$\gather
\Delta (k^{\pm 1}) = k^{\pm 1}\otimes k^{\pm 1}, 
\Delta (e) = e\otimes 1 + k\otimes e,
\Delta (f) = f\otimes k^{-1} + 1\otimes f; \tag 2.2\\
S(k^{\pm 1}) = k^{\mp 1}, \quad S(e)=-k^{-1}e, \quad S(f) = -fk; \tag 2.3\\
\varepsilon (k^{\pm 1})=1, \quad \varepsilon (e)=\varepsilon (f) = 0. \tag 2.4
\endgather
$$
The quantum supergroup $\Cal{A}=\Bbb{C}_{q}[osp(1,2)]$ is the $\Bbb{Z}_{2}$-graded associative algebra over $\Bbb{C}$ defined by generators $a, b, c, d, \sigma$ with relations
$$
\gathered
ab = tba, \qquad ac = tca, \qquad bc = -cb, \\
bd = -tdb, \qquad cd = -tdc, \qquad ad - da = (t^{-1}-t)bc,\\
ad + tbc = \sigma, \quad \sigma^{2}=1. 
\endgathered \tag 2.5
$$
The $\Bbb{Z}_{2}$-grading on $\Cal{A}$ is given by $p(a)=p(d)=p(\sigma)=0$ and $p(b)=p(c)=1$.  The $\Bbb{Z}_{2}$-graded Hopf algebra structure on $\Cal{A}$ is defined by
$$
\gathered
\Delta \left( \matrix a & b & 0\\ c & d & 0\\ 0 & 0 & \sigma \endmatrix \right)=
\left( \matrix a & b & 0\\ c & d & 0\\ 0 & 0 & \sigma \endmatrix \right)\otimes
\left(\matrix a & b & 0\\ c & d & 0\\ 0 & 0 & \sigma \endmatrix \right),\\
S\left( \matrix a & b & 0\\ c & d & 0\\ 0 & 0 & \sigma \endmatrix \right)=
\left( \matrix d\sigma & -t^{-1}b\sigma & 0\\ tc\sigma & a\sigma & 0\\ 0 & 0 & \sigma \endmatrix \right),\\ 
\varepsilon \left(\matrix a & b & 0\\ c & d & 0\\ 0 & 0 & \sigma \endmatrix \right) = \left(\matrix 1 & 0 & 0\\ 0 & 1 & 0\\ 0 & 0 & 1 \endmatrix \right),
\endgathered \tag 2.6
$$
\par
The dual pairing between $\Cal{U}$ and $\Cal{A}$ is provided by
$$
\gathered
k^{\pm 1}\left( \matrix a & b & 0\\ c & d & 0\\ 0 & 0 & \sigma \endmatrix \right)=
\left( \matrix t^{\pm 1} & 0 & 0\\ 0 & -t^{\mp 1} & 0\\ 0 & 0 & -1 \endmatrix \right),\\
e\left( \matrix  a & b & 0 \\ c & d & 0\\ 0 & 0 & \sigma \endmatrix \right)=
\left( \matrix 0 & \frac{t-t^{-1}}{q-q^{-1}} & 0\\ 0 & 0 &0\\0 & 0 & 0 \endmatrix \right), \quad
f\left( \matrix  a & b & 0\\ c & d & 0\\ 0 & 0 & \sigma \endmatrix \right)=
\left( \matrix 0 & 0 &0\\ 1 & 0 & 0\\ 0& 0 & 0 \endmatrix \right),
\endgathered\tag 2.7
$$
together with the requirements that $k^{\pm 1}$ are algebra homomorphisms and 
$$
\gathered
e(xy)= e(x)\varepsilon (y) + (-1)^{p(x)}k(x)e(y), \qquad e(1)=0,\\
f(xy)= f(x)k^{-1}(y) + (-1)^{p(x)}\varepsilon(x)f(y), \qquad f(1)=0,
\endgathered\tag 2.8
$$
for all $x,y\in \Cal{A}$. For $u\in\Cal{U}$, we also have
$$
\Delta (u)(x\otimes y) = u(xy), \quad \varepsilon(u) = u(1), \quad S(u)(x)=u(S(x)). \tag 2.9
$$
\par
Recall the following definition of the Gauss' binomial coefficients:
$$
(u;v)_{m}=\prod_{k=0}^{m-1}(1-uv^{k}), \quad 
\left(\matrix m\\ n \endmatrix\right)_{v}=\frac{(v;v)_{m}}{(v;v)_{n}(v;v)_{m-n}}. 
$$
\par
For each $\ell\in\frac{1}{2}\Bbb{Z}_{+}$, let $I_{\ell} =\{-\ell,-\ell + 1, ...,\ell\}$ and define two $\Bbb C$-vector subspaces $V^{L}_{\ell}$ and $V^{R}_{\ell}$ of $\Cal{A}$ by ([Z, Section 5])
$$
V^{L}_{\ell}=\bigoplus_{i\in I_{\ell}}\Bbb{C}\xi^{(\ell)}_{i}\qquad\text{and}\qquad
V^{R}_{\ell}=\bigoplus_{i\in I_{\ell}}\Bbb{C}\eta^{(\ell)}_{i},\tag 2.10
$$
where 
$$
\xi^{(\ell)}_{i} =\bold{i}^{[\frac{\ell+i}{2}]}\left(\matrix 2\ell\\ \ell+i\endmatrix\right)^{\frac{1}{2}}_{t^{-2}}
a^{\ell -i}c^{\ell +i},\quad\text{and}\quad
\eta^{(\ell)}_{i} =\bold{i}^{[\frac{\ell+i}{2}]}\left(\matrix 2\ell\\ \ell+i\endmatrix\right)^{\frac{1}{2}}_{t^{-2}}
a^{\ell -i}b^{\ell +i}.\tag 2.11
$$
We shall denote the subspace spanned by the elements $\xi^{(\ell)}_{i}\sigma$ (resp. $\eta^{(\ell)}_{i}\sigma$) ($i\in I_{\ell}$) by $V^{L}_{\ell}\sigma$ (resp. $V^{R}_{\ell}\sigma$).  Then each $V^{L}_{\ell}\sigma^{s}$ (resp. $V^{R}_{\ell}\sigma^{s}$), $s=0,1$, forms an irreducible left (resp. right) $\Cal{A}$-subcomodule of $\Cal{A}$.  Let $m^{(\ell)}_{ij}\sigma^{s}$ ($i,j\in I_{\ell}$) be the associated matrix elements. Then we have 
$$
\gathered
\Delta(\xi^{(\ell)}_{i}\sigma^{s})=\sum_{j\in I_{\ell}}m^{(\ell)}_{ij}\sigma^{s}\otimes\xi^{(\ell)}_{j}\sigma^{s}, \qquad i\in I_{\ell},\\
\Delta(\eta^{(\ell)}_{i}\sigma^{s})=\sum_{j\in I_{\ell}}\eta^{(\ell)}_{j}\sigma^{s}\otimes m^{(\ell)}_{ji}\sigma^{s}, \\
\Delta(m^{(\ell)}_{ij}\sigma^{s})=\sum_{k\in I_{\ell}}m^{(\ell)}_{ik}\sigma^{s}\otimes m^{(\ell)}_{kj}\sigma^{s}, \qquad \text{and}\qquad
\varepsilon(m^{(\ell)}_{ij}\sigma^{s})=\delta_{ij}.
\endgathered
\tag 2.12
$$
\par
Recall the definition of the little $q$-Jacobi polynomials (see [MNU, 2.2]):
$$
P^{(\alpha,\beta)}_{n}(z;q)=\sum_{r\ge 0}\frac{(q^{-n};q)_{r}(q^{\alpha+\beta+n+1};q)_{r}}{(q;q)_{r}(q^{\alpha+1};q)_{r}}(qz)^{r}, 
$$
where $\alpha, \beta \in \Bbb{Z}, n\in \Bbb{Z}_{+}$.  Let
$$
N^{(\ell)}_{ij} = \bold{i}^{[\frac{\ell+i}{2}]-[\frac{\ell+j}{2}]} t^{(\ell+j)(i-j)}\left(\matrix \ell + i\\i-j\endmatrix\right)^{\frac{1}{2}}_{t^{-2}}\left(\matrix \ell -j\\i-j\endmatrix\right)^{\frac{1}{2}}_{t^{-2}}.
$$
One can compute the matrix elements by the following theorem ([Z, Thm. 6.1]):
\par
\proclaim{Theorem 2.1}
Let $\zeta = tbc\sigma$. For $\ell\in\frac{1}{2}\Bbb{Z}_{+}$, $i,j\in I_{\ell}$, we have
\par
(1) If $i+j \le 0$, $i\ge j$, then 
$$
m^{(\ell)}_{ij}= a^{-i-j}c^{i-j}\sigma^{\ell+j}N^{(\ell)}_{ij}P^{(i-j,-i-j)}_{\ell+j}(\zeta; t^{-2}).\tag 2.13
$$
\par
(2) If $i+j \le 0$, $j\ge i$, then
$$
m^{(\ell)}_{ij}= (-1)^{[\frac{j-i}{2}]} 
a^{-i-j}b^{j-i}\sigma^{\ell+i}N^{(\ell)}_{ji}P^{(j-i,-i-j)}_{\ell+i}(\zeta; t^{-2}).\tag 2.14
$$
\par
(3) If $i+j \ge 0$, $j\ge i$, then
$$
m^{(\ell)}_{ij}= (-1)^{[\frac{j-i}{2}]} 
N^{(\ell)}_{-i,-j}P^{(j-i,i+j)}_{\ell-j}(\zeta; t^{-2})
b^{j-i}d^{i+j}\sigma^{\ell-j}.\tag 2.15
$$
\par
(4) If $i+j \ge 0$, $i\ge j$, then
$$
m^{(\ell)}_{ij}= 
N^{(\ell)}_{-j,-i}P^{(i-j,i+j)}_{\ell-i}(\zeta; t^{-2})
c^{i-j}d^{i+j}\sigma^{\ell-i}.\tag 2.16
$$
\endproclaim
\par
For $\ell\in I_{\ell}$, $s=0,1$, let $M_{\ell}\sigma^{s}=\sum_{i,j\in I_{\ell}}\Bbb{C}m^{(\ell)}_{i,j}\sigma^{s}$, then $M_{\ell}\sigma^{s}$ is a bi-subcomodule of $\Cal{A}$, and we can define an $\Cal{A}$-bicomodule isomorphism $\Phi_{\ell} : V^{L}_{\ell}\sigma^{s}\otimes V^{R}_{\ell}\sigma^{s} \longrightarrow M_{\ell}\sigma^{s}$ by
$$
\Phi_{\ell}(\xi^{(\ell)}_{i}\sigma^{s}\otimes\eta^{(\ell)}_{j}\sigma^{s})= m^{(\ell)}_{ij}\sigma^{s}, \qquad
i,j\in I_{\ell}. \tag 2.17
$$
\par
\proclaim{Theorem 2.2} ([Z], Theorem 1.1 and Theorem 5.1) (1) The set of monomials $a^{i}b^{j}c^{k}d^{l}\sigma^{s}$, $i,j,k,l \in \Bbb{Z}_{+}$, $i=0$ or $l=0$, and $s=0,1$, form a basis of $\Cal{A}$ over $\Bbb{C}$. 
\par
(2) Any finite dimensional irreducible left (or right) $\Cal{A}$-comodule is isomorphic to $V^{L}_{\ell}\sigma^{s}$ (or $V^{R}_{\ell}\sigma^{s}$) ($s=0,1$) for some $\ell\in \frac{1}{2}\Bbb{Z}_{+}$, and $\Cal{A}$ is decomposed into the direct sum of $\Cal{A}$-bicomodules
$$
\Cal{A} =\bigoplus_{\ell\in\frac{1}{2}\Bbb{Z}_{+}}(M_{\ell}\oplus M_{\ell}\sigma).\tag 2.18
$$
\endproclaim
\par
\remark{Remark}We can define two hermitian forms $<,>_{R}$ and $<,>_{L}$ on $\Cal{A}$ such that the decomposition (2.18) is an orthogonal decomposition with respect to both forms, for details see [Z, Section 7].
\endremark
\head
3. Primitive elements and homogeneous spaces
\endhead
\par
Recall that for a pair of group-like elements $(g, g')$ in a Hopf algebra $H$, an element $x\in H$ is called $(g, g')$-primitive if $\Delta(x) = x\otimes g + g'\otimes x$.  If $x$ is $(g, g')$-primitive, then the one-dimensional subspace spanned by $x$ is a coideal of $H$.  We will study the homogeneous superspaces of $\Cal{A}$ defined by the $(g, g')$-primitive elements of $\Cal{U}$, so we need to describe the $(g,g')$-primitive elements.
\proclaim{Proposition 3.1} 
The group-like elements of $\Cal{U}$ are $k^{m}, m\in \Bbb{Z}$. For $m, n\in \Bbb{Z}$, denote by $P(m,n)$ the subspace of all $(k^{m}, k^{n})$-primitive elements of $\Cal{U}$.
\par
(1) If $m=n$, then $P(m,n) = 0$.
\par
(2) If $n - m \ne 0, 1$, then $P(m,n) = \Bbb{C}(k^{m} - k^{n})$.
\par
(3) $P(0,1) = \Bbb{C}e \oplus \Bbb{C}fk \oplus \Bbb{C}(k-1)$.  
\par
(4) If $n - m = 1$, then $P(m,n) = k^{m}P(0,1)$.
\endproclaim  
\demo{Proof}
The monomials $e^{r}f^{s}k^{t}$ ($r, s \in \Bbb{Z}_{+}, t\in \Bbb{Z}$) form a basis of $\Cal{U}$, we consider the action of $\Delta$ on these monomials.  By using the formula (which holds in any associative algebra for elements $x$ and $y$ satisfying $xy = vyx$)
$$
(x+y)^{m}=\sum_{k=0}^{m}\left(\matrix m\\ k \endmatrix\right)_{v^{-1}}x^{k}y^{m-k},
$$
one can prove
$$
\Delta(e^{r}f^{s}k^{t}) =\sum_{i=0}^{r}\sum_{j=0}^{s}c(r,s,i,j;q) e^{i}f^{s-j}k^{r-i+t}\otimes e^{r-i}f^{j}k^{j-s+t}, \tag 3.1
$$
where
$$
c(r,s,i,j;q)= (-1)^{(r-i+j)(s-j)}q^{(i-r)(s-j)} \left(\matrix r\\ i \endmatrix\right)_{-q}\left(\matrix s\\ j \endmatrix\right)_{-q}.
$$
Order the monomials $e^{r}f^{s}k^{t}$ lexicographically, and express a group-like element $u\in \Cal{U}$ as a linear combination of this basis.  If $ce^{r}f^{s}k^{t}$ is the leading term of $u$, then using (3.1) one can see that $c^{2}e^{r}f^{s}k^{t}\otimes e^{r}f^{s}k^{t}$ shows up in $\Delta(u)$ only if $r = s = 0$, i.e. $u$ must be a linear combination of the $k^{t}$.  Then it follows that the group-like elements are $k^{t}$ ($t\in \Bbb{Z}$).
\par
Suppose that $u$ is a $(k^{m}, k^{n})$-primitive element. If both $r, s > 0$ in the leading term $e^{r}f^{s}k^{t}$ of $u$, then from
$$
\Delta (e^{r}f^{s}k^{t}) = e^{r}f^{s}k^{t}\otimes k^{-s+t} + (-1)^{s-1}\left(\matrix s\\ 1 \endmatrix\right)_{-q} e^{r}f^{s-1}k^{t} \otimes fk^{1-s+t} + \cdots, 
$$
one can see that the second term can not be cancelled in $\Delta(u)$, so $s = 0$ or $r = 0$.  Similar arguments further show that $u$ must be a linear combination of $ek^{t}, fk^{t'}, k^{j}$ ($t, t', j\in \Bbb{Z}$).  Now the statements in the proposition follow easily by considering the action of $\Delta$ on a linear combination of these elements. 
\enddemo
\par
For each $\ell \in \frac{1}{2}\Bbb{Z}_{+}$, $\Cal{U}$ has up to equivalence a unique representation of dimension $2\ell + 1$. The finite-dimensional representations of $\Cal{U}$ can be realized by the subspaces $V^{L}_{\ell}\sigma^{s}$ (as right $\Cal{U}$-modules) or $V^{R}_{\ell}\sigma^{s}$ (as left $\Cal{U}$-modules) of $\Cal{A}$.  The left actions of $e, f, k^{\pm 1}$ on the basis $\eta^{(\ell)}_{j}\sigma^{s}$ can be computed by using the formulas in Theorem 2.1, and we have
$$ 
\aligned
L_{k^{\pm 1}}(\eta^{(\ell)}_{j}\sigma^{s}) &= (-1)^{\ell + j+s}t^{\mp 2j}\eta^{(\ell)}_{j}\sigma^{s},\\
L_{e}(\eta^{(\ell)}_{j}\sigma^{s}) &= (-1)^{\ell - j + 1}\bold{i}^{[\frac{\ell +j-1}{2}]-[\frac{\ell + j}{2}]}t^{\ell - j + 1}\\
&\times \frac{(1-t^{-2(\ell + j)})^{\frac{1}{2}}(1-t^{-2(\ell - j +1)})^{\frac{1}{2}}}{q - q^{-1}}\eta^{(\ell)}_{j-1}\sigma^{s},\\ 
L_{f}(\eta^{(\ell)}_{j}\sigma^{s}) &= (-1)^{1+s}\bold{i}^{[\frac{\ell +j+1}{2}]-[\frac{\ell + j}{2}]}t^{\ell + j }\\
&\times\frac{(1-t^{-2(\ell - j)})^{\frac{1}{2}}(1-t^{-2(\ell + j +1)})^{\frac{1}{2}}}{1-t^{-2}}\eta^{(\ell)}_{j+1}\sigma^{s},\endaligned
\tag 3.2
$$
with the convention that $\eta^{(\ell)}_{-\ell - 1}$ and $\eta^{(\ell)}_{\ell + 1}$ are zero.  The right actions on the basis $\xi^{(\ell)}_{i}\sigma^{s}$ are given by
$$\aligned
R_{k^{\pm 1}}(\xi^{(\ell)}_{i}\sigma^{s}) &= (-1)^{\ell + i+s}t^{\mp 2i}\xi^{(\ell)}_{i}\sigma^{s},\\
R_{e}(\xi^{(\ell)}_{i}\sigma^{s}) &= \bold{i}^{[\frac{\ell + i}{2}]-[\frac{\ell + i+1}{2}]}t^{\ell - i}\frac{(1-t^{-2(\ell -i)})^{\frac{1}{2}}(1-t^{-2(\ell +i +1)})^{\frac{1}{2}}}{q - q^{-1}}\xi^{(\ell)}_{i+1}\sigma^{s},\\
R_{f}(\xi^{(\ell)}_{i}\sigma^{s}) &= (-1)^{\ell + i+s-1}\bold{i}^{[\frac{\ell +i}{2}]-[\frac{\ell + i-1}{2}]}t^{\ell + i-1}\\
&\times \frac{(1-t^{-2(\ell +i)})^{\frac{1}{2}}(1-t^{-2(\ell -i +1)})^{\frac{1}{2}}}{1-t^{-2}}\xi^{(\ell)}_{i-1}\sigma^{s}.\endaligned
\tag 3.3
$$
\par
According to Proposition 3.1, the basic $(1,k)$-primitive elements are linear combinations of $e$, $fk$ and $k - 1$.  Let 
$$
x_{\alpha\beta\gamma} = \alpha e + \beta fk + \gamma (k-1), \tag 3.4
$$
where $\alpha, \beta, \gamma \in \Bbb{C}$, not all zero. If one of $\alpha, \beta, \gamma$ is zero, we shall omit it from the index, i.e. $x_{\alpha\beta}=\alpha e + \beta fk$ and so on. We consider the invariant subspace of $\Cal{A}$ defined by $x_{\alpha\beta\gamma}$. Note that $\varepsilon (x_{\alpha\beta\gamma})=0$, so the invariant subspace is the null space of $R_{x_{\alpha\beta\gamma}}$ (or $L_{x_{\alpha\beta\gamma}}$).
\par
 By (3.3), for $\ell = \frac{1}{2}$, the right action of $x_{\alpha\beta\gamma}$ on the basis $\xi^{(\ell)}_{i}\sigma^{s}$ corresponds to the matrix
$$
\left( \matrix \gamma ((-1)^{s}t-1) & -\beta t^{-1} \\ \alpha \frac{t - t^{-1}}{q - q^{-1}} & \gamma ((-1)^{1+s}t^{-1}-1)\endmatrix \right).
$$
The determinant of this matrix is equal to 
$$
(-1)^{s}\gamma^{2}(-t + t^{-1}) + \alpha\beta t^{-1}\frac {t - t^{-1}}{q - q^{-1}}.
$$
\par
If $\gamma = 0$, then $\alpha\beta$ must be $0$ in order for $x_{\alpha\beta\gamma}$ to have nontrivial invariant subspace in $V^{L}_{\frac{1}{2}}\sigma^{s}$, i.e. $x_{\alpha\beta\gamma}$ is a nonzero multiple of $e$ or $fk$.  Note that in both cases, there is a one-dimensional invariant subspace for every $\ell \in \frac{1}{2}\Bbb{Z}_{+}$, given by $\xi^{(\ell)}_{\ell}\sigma^{s}$ or $\xi^{(\ell)}_{-\ell}\sigma^{s}$ (respetively).
\par
If $\gamma \ne 0$, we can assume $\gamma = 1$ and chose $\alpha$, $\beta$ so that
$$
\alpha\beta = (-1)^{s}t(q - q^{-1}).\tag 3.5
$$ 
Then $x_{\alpha\beta\gamma}$ has a one-dimensional invariant subspace in $V^{L}_{\frac{1}{2}}\sigma^{s}$, spanned by the element
$$
((-1)^{s}t^{-1} + 1)\xi^{(\frac{1}{2})}_{-\frac{1}{2}} + \alpha\frac{t - t^{-1}}{q - q^{-1}}\xi^{(\frac{1}{2})}_{\frac{1}{2}}. \tag 3.6
$$
\par
However if $\gamma \ne 0$, then there is no nontrivial invariants for higher dimensions in general.  Consider the case $\ell = 1$. The matrix of the right action of $x_{\alpha\beta\gamma}$ on $V^{L}_{1}$ is
$$
R_{x_{\alpha\beta\gamma}}(3) = \left( \matrix -\gamma (q+1) & -\beta (1-q^{-1})^{\frac{1}{2}} & 0 \\ -\alpha (1-q^{-1})^{-\frac{1}{2}} & -2\gamma & -\bold{i}\beta t^{-1}(1-q^{-1})^{\frac{1}{2}}\\ 0 &  \bold{i}\alpha t^{-1}(1-q^{-1})^{-\frac{1}{2}} & -\gamma (1+q^{-1})\endmatrix \right), \tag 3.7
$$
with the determinant 
$$
\det(R_{x_{\alpha\beta\gamma}}(3)) = -2\frac{\gamma^{3}(1+q)^{2}}{q},
$$
and $\det(R_{x_{\alpha\beta\gamma}}(3)) = 0$ if and only if $\gamma = 0$ (under our assumption on $q$).  The matrix of the right action of $x_{\alpha\beta\gamma}$ on $V^{L}_{1}\sigma$ is
$$
R_{x_{\alpha\beta\gamma}}(3,\sigma) = \left( \matrix \gamma (q-1) & -\beta (1-q^{-1})^{\frac{1}{2}} & 0 \\ -\alpha (1-q^{-1})^{-\frac{1}{2}} & 0 & -\bold{i}\beta t^{-1}(1-q^{-1})^{\frac{1}{2}}\\ 0 &  \bold{i}\alpha t^{-1}(1-q^{-1})^{-\frac{1}{2}} & \gamma (q^{-1}-1)\endmatrix \right), \tag 3.8
$$
with the determinant
$$
\det(R_{x_{\alpha\beta\gamma}}(3,\sigma)) = 2\frac{\alpha\beta\gamma(q-1)}{q},
$$
and $\det R_{x_{\alpha\beta\gamma}}(3,\sigma) = 0$ if and only if $\alpha\beta\gamma =0$.  Further computations show that under condition (3.5), the parameter $q$ must be a root of a polynomial (depending on $\ell$) in order for the determinant to be 0.  The left action can be discussed similarly.  We summarize the above discussion in the following theorem.
\proclaim{Theorem 3.2}
Let $x_{\alpha\beta\gamma}$ be defined by (3.4).  Then in order for the right (or left) action of $x_{\alpha\beta\gamma}$ on $\Cal{A}$ to have nontrivial invariants for all $\ell\in \frac{1}{2}\Bbb{Z}_{+}$, $x_{\alpha\beta\gamma}$ must be scalar multiples of $e$ or $fk$.  The right invariants of the coideal $\Bbb{C}e$ are the linear combinations of $m^{(\ell)}_{\ell j}\sigma^{s}$ ($\ell \in \frac{1}{2}\Bbb{Z}_{+}; i,j\in I_{\ell}$), and the right invariants of the coideal $\Bbb{C}fk$ are the linear combinations of $m^{(\ell)}_{-\ell j}\sigma^{s}$ ($\ell \in \frac{1}{2}\Bbb{Z}_{+}; i,j\in I_{\ell}$).
\endproclaim
\demo{Proof}
Use the $\Cal A$-bimodule isomorphism $\Phi_{\ell}$ defined by (2.17).
\enddemo
\par
If $\gamma=0$ and $\alpha\beta\ne 0$, then $R_{x_{\alpha\beta}}$ has a one-dimensional invariant subspace in $V^{L}_{1}\sigma^{s}$ ($s=0,1$). To treat the general cases, we need the following lemma, which can be proved by standard linear algebra arguments.
\proclaim{Lemma 3.3} Let
$$
A = \pmatrix 
\format \c&\quad \c&\quad \c&\quad \c&\quad \c&\c&\, \c\\
0&\alpha_{1}&0&&&\hdots&0\\
\beta_{1}&0&\alpha_{2}&0&&\hdots&0\\
0&\beta_{2}&0&\alpha_{3}&0&\hdots&0\\
\vdots&&&&\ddots&&\vdots\\
0&&\hdots&&0&&\alpha_{n-1}\\
0&&\hdots&0&\beta_{n-1}&&0
\endpmatrix.
$$
Then
$$
\det(A) = \cases 0, &\text{$n$ odd,}\\ 
(-1)^{\frac{n}{2}}\alpha_{1}\beta_{1}\alpha_{3}\beta_{3}\cdots
\alpha_{n-1}\beta_{n-1}, &\text{$n$ even.}\endcases
$$
Hence if all $\alpha_{i}$ and $\beta_{i}$ are nonzero, then the linear system 
$$AX=0, \qquad X=(x_{1},...,x_{n})^{T},
$$
has nontrivial solution if and only if $n$ is odd, and if this is the case, the solutions are scalar multiples of the column vector $X$ with components 
$$
x_{i} = \cases 0, &\text{$i$ even,}\\
(-1)^{\frac{i-1}{2}}\frac{\prod_{j=1}^{\frac{i-1}{2}}\beta_{2j-1}}
{\prod_{j=1}^{\frac{i-1}{2}}\alpha_{2j}}, &\text{$i$ odd,}
\endcases
$$ 
with the convention that $x_{1} = 1$. 
\endproclaim  
\par
For $i\in I_{\ell}$, we define $a_{i}, b_{i}$, and $c_{i}$ by 
$$
R_{k-1}(\xi^{(\ell)}_{i}\sigma^{s}) = a_{i}\xi^{(\ell)}_{i}\sigma^{s},\quad
R_{e}(\xi^{(\ell)}_{i}\sigma^{s}) = b_{i}\xi^{(\ell)}_{i+1}\sigma^{s},\quad
R_{fk}(\xi^{(\ell)}_{i}\sigma^{s}) = c_{i}\xi^{(\ell)}_{i-1}\sigma^{s}. \tag 3.9
$$
According to (3.3), $b_{i}\ne 0$ for all $i < \ell$ and $c_{i}\ne 0$ for all $-\ell < i$.  Since the matrix of $R_{x_{\alpha\beta}}$ on the basis $\xi^{(\ell)}_{-\ell}\sigma^{s}, \xi^{(\ell)}_{-\ell +1}\sigma^{s},..., \xi^{(\ell)}_{\ell}\sigma^{s}$ of $V^{L}_{\ell}\sigma^{s}$ is
$$
R_{x_{\alpha\beta}}(\ell,\sigma^{s}) = \pmatrix 
\format \c&\quad \c&\quad \c&\quad \c&\quad \c&\c&\, \c\\
0&\beta c_{-\ell+1}&0&&&\hdots&0\\
\alpha b_{-\ell}&0&\beta c_{-\ell+2}&0&&\hdots&0\\
0&\alpha b_{-\ell+1}&0&\beta c_{-\ell+3}&0&\hdots&0\\
\vdots&&&&\ddots&&\vdots\\
0&&\hdots&&0&&\beta c_{\ell}\\
0&&\hdots&0&\alpha b_{\ell-1}&&0
\endpmatrix,
$$ 
by Lemma 3.3, $x_{\alpha\beta}$ has nontrivial invariant in $V^{L}_{\ell}\sigma^{s}$ if and only if $2\ell + 1$ is odd, that is, if and only if $\ell\in \Bbb{Z}_{+}$, and the invariants are scalar multiples of 
$$
\rho^{(\ell)}_{\alpha\beta} =\sum_{i\in I_{\ell}}x_{i}\xi^{(\ell)}_{i}\sigma^{s},\tag 3.10
$$
where 
$$
x_{i} = \cases 
1, & \text{if $i=-\ell$},\\
0, &\text{if $i=-\ell+2k+1$, \quad $0\le k\le \ell-1$},\\
(-\alpha\beta^{-1})^{k}\frac{\prod_{j=0}^{k-1}b_{-\ell+2j}}{\prod_{j=1}^{k}c_{-\ell+2j}}, &\text{if $i=-\ell+2k$, \quad $1\le k\le \ell$}.
\endcases\tag 3.11
$$
\par
Note that
$$
\gathered
\prod_{j=0}^{k-1}b_{-\ell+2j} = (-q)^{k\ell-\frac{k(k-1)}{2}}\frac{(q^{-2\ell};q^{2})^{\frac{1}{2}}_{k}(-q^{-1};q^{-2})^{\frac{1}{2}}_{k}}{(q-q^{-1})^{k}},\\
\prod_{j=1}^{k}c_{-\ell+2j} = (-\bold{i})^{k}(-q)^{k\ell-\frac{k(k+1)}{2}}\frac{(-q^{-2\ell+1};q^{2})^{\frac{1}{2}}_{k}(q^{-2};q^{-2})^{\frac{1}{2}}_{k}}{(t-t^{-1})^{k}},
\endgathered\tag 3.12
$$ 
so we can write
$$\split
\frac{\prod_{j=0}^{k-1}b_{-\ell+2j}}{\prod_{j=1}^{k}c_{-\ell+2j}}&= \left(\frac{\bold{i}q}{t+t^{-1}}\right)^{k}\frac{(q^{-2\ell};q^{2})^{\frac{1}{2}}_{k}(-q^{-1};q^{-2})^{\frac{1}{2}}_{k}}{(-q^{-2\ell+1};q^{2})^{\frac{1}{2}}_{k}(q^{-2};q^{-2})^{\frac{1}{2}}_{k}}\\
&= \bold{i}^{k^{2}}\left(\frac{q^{2}}{1-q}\right)^{k}\frac{(q^{-2\ell};q^{2})^{\frac{1}{2}}_{k}(-q;q^{2})^{\frac{1}{2}}_{k}}
{(-q^{-2\ell+1};q^{2})^{\frac{1}{2}}_{k}(q^{2};q^{2})^{\frac{1}{2}}_{k}}.
\endsplit\tag 3.13
$$
\par
Now assume $\gamma\ne 0$ but $\alpha\beta = 0$.  For the case $\alpha = 0$, the matrix of $x_{\beta\gamma}$ on $V^{L}_{\ell}\sigma^{s}$ is
$$
R_{x_{\beta\gamma}}(\ell, s) = \pmatrix 
\format \c&\quad \c&\quad \c&\quad \c&\quad \c&\c&\, \c\\
\gamma a_{-\ell}&\beta c_{-\ell+1}&0&&&\hdots&0\\
0&\gamma a_{-\ell+1}&\beta c_{-\ell+2}&0&&\hdots&0\\
0&0&\gamma a_{-\ell+2}&\beta c_{-\ell+3}&0&\hdots&0\\
\vdots&&&&\ddots&&\vdots\\
0&&\hdots&&\gamma a_{\ell-1}&&\beta c_{\ell}\\
0&&\hdots&0&0&&\gamma a_{\ell}
\endpmatrix,
$$
and $x_{\beta\gamma}$ has nontrivial invariant in $V^{L}_{\ell}\sigma^{s}$ $\Leftrightarrow\prod_{i\in I_{\ell}}a_{i}=0 \Leftrightarrow 0\in I_{\ell}$ and $\ell + s \in 2\Bbb{Z}$ (by (3.3)).  This last condition is equivalent to $\ell = 2r+s, r \in \Bbb{Z}_{+}$.  Under this condition, the one-dimensional invariant subspace is spanned by
$$
\rho^{(\ell)}_{\beta\gamma} =\sum_{i\in I_{\ell}}x_{i}\xi^{(\ell)}_{i}\sigma^{s},\tag 3.14
$$
where 
$$
x_{i} = \cases 
1, & \text{if $i=0$},\\
0, &\text{if $1\le i\le \ell$},\\
(-\beta\gamma^{-1})^{-i}\frac{\prod_{j=i+1}^{0}c_{j}}{\prod_{j=i}^{-1}a_{j}}, &\text{if $-\ell\le i\le -1$}.
\endcases\tag 3.15
$$
Note that we can write
$$
\frac{\prod_{j=i+1}^{0}c_{j}}{\prod_{j=i}^{-1}a_{j}} = \bold{i}^{[\frac{\ell}{2}]-[\frac{\ell+i}{2}]+i}\frac{((-q)^{-\ell};-q)^{\frac{1}{2}}_{-i}((-q)^{\ell+1};-q)^{\frac{1}{2}}_{-i}}{(1+q)^{-i}(q;q)_{-i}}.
\tag 3.16
$$
\par
Similarly, $x_{\alpha\gamma}$ has nontrivial invariant subspace if and only if $\ell = 2r+s, r \in \Bbb{Z}_{+}$, and if this is the case, the one-dimensional invariant subspace is spanned by
$$
\rho^{(\ell)}_{\alpha\gamma} =\sum_{i\in I_{\ell}}x'_{i}\xi^{(\ell)}_{i}\sigma^{s},\tag 3.17
$$
where 
$$
x'_{i} = \cases 
0, &\text{if $-\ell\le i\le -1$},\\
1, & \text{if $i=0$},\\
(-\alpha\gamma^{-1})^{i}\frac{\prod_{j=0}^{i-1}b_{j}}{\prod_{j=1}^{i}a_{j}}, &\text{if $1\le i\le \ell$}.
\endcases\tag 3.18
$$
We also have
$$
\frac{\prod_{j=0}^{i-1}b_{j}}{\prod_{j=1}^{i}a_{j}} = \bold{i}^{[\frac{\ell}{2}]-[\frac{\ell+i}{2}]+i}t^{i}
\frac{((-q)^{-\ell};-q)^{\frac{1}{2}}_{i}((-q)^{\ell+1};-q)^{\frac{1}{2}}_{i}}{(q-q^{-1})^{i}(q;q)_{i}}.
\tag 3.19
$$
\par
We now summarize our discussion as follows.
\proclaim{Proposition 3.4}
The right actions of the elements $x_{\alpha\beta}$, $x_{\alpha\gamma}$, and $x_{\beta\gamma}$ have nontrivial invariants in $V^{L}_{\ell}\sigma^{s}$ only if $\ell\in \Bbb{Z}_{+}$.  If $\ell\in \Bbb{Z}_{+}$, then $x_{\alpha\beta}$ has a one-dimensional invariant subspace in $V^{L}_{\ell}\sigma^{s}$ for both $s=0,1$ with a basis given by (3.10) and (3.11), while $x_{\alpha\gamma}$ (resp. $x_{\beta\gamma}$) has nontrivial invariants only if $\ell+s$ is even, and if this is the case, the one-dimensional invariant subspace is given by (3.14) and (3.15) (resp. (3.17) and (3.18)).
\endproclaim
\par
From this proposition, we immediately obtain the following theorem.
\proclaim{Theorem 3.5} Let $\alpha, \beta$, and $\gamma$ be nonzero complex numbers.
\par
(1) The right invariant subalgebra of $\Cal{A}$ corresponding to $x_{\alpha\beta}$ is spanned by 
$$
\sum_{i\in I_{\ell}}x_{i}m^{(\ell)}_{ij}\sigma^{s},\quad \ell\in \Bbb{Z}_{+},\quad j\in I_{\ell},\quad s=0,1,\tag 3.20
$$
where $x_{i}$ are defined by (3.11).
\par
(2) The right invariant subalgebra of $\Cal{A}$ corresponding to $x_{\alpha\gamma}$ (resp. $x_{\beta\gamma}$) is spanned by
$$
\sum_{i\in I_{\ell}}x_{i}m^{(\ell)}_{ij}\sigma^{s},\quad j\in I_{\ell},\tag 3.21
$$
where $\ell\in\Bbb{Z}_{+}$ and $s=0,1$ satisfy $\ell+s \in 2\Bbb{Z}_{+}$, and $x_{i}$ are defined by (3.15) (resp. (3.18)).
\endproclaim
\par
The left action can be discussed similarly by using (3.2), we give the corresponding results for the left action here and omit the details. 
\par
\proclaim{Proposition 3.6}
The left actions of elements $x_{\alpha\beta}$, $x_{\alpha\gamma}$, and $x_{\beta\gamma}$ have nontrivial invariants in $V^{R}_{\ell}\sigma^{s}$ only if $\ell\in \Bbb{Z}_{+}$.  Assume $\ell\in \Bbb{Z}_{+}$, let
$$
\lambda^{(\ell)}_{\alpha\beta\gamma} =\sum_{j\in I_{\ell}}y_{j}\eta^{(\ell)}_{j}\sigma^{s}.\tag 3.22
$$
\par
(1) The element $x_{\alpha\beta}$ has a one-dimensional invariant subspace in $V^{R}_{\ell}\sigma^{s}$ for both $s=0,1$ spanned by the element (3.22) with
$$
y_{j} = \cases 
1, & \text{if $j=-\ell$},\\
0, &\text{if $j=-\ell+2k+1$, \quad $0\le k\le \ell-1$},\\
(-\alpha^{-1}\beta)^{k}\frac{\prod_{i=0}^{k-1}c'_{-\ell+2i}}{\prod_{i=1}^{k}b'_{-\ell+2i}}, &\text{if $j=-\ell+2k$, \quad $1\le k\le \ell$},
\endcases\tag 3.23
$$
where
$$
\frac{\prod_{i=0}^{k-1}c'_{-\ell+2i}}{\prod_{i=1}^{k}b'_{-\ell+2i}}
= \bold{i}^{k}q^{k}(1-q)^{k}\frac{(q^{-2\ell};q^{2})^{\frac{1}{2}}_{k}(-q;q^{2})^{\frac{1}{2}}_{k}}
{(-q^{-2\ell+1};q^{2})^{\frac{1}{2}}_{k}(q^{2};q^{2})^{\frac{1}{2}}_{k}}.
\tag 3.24
$$
\par
(2) The elements $x_{\beta\gamma}$ and $x_{\alpha\gamma}$ have nontrivial invariants only if $\ell+s$ is even, and if this is the case, the one-dimensional invariant subspace is given by (3.22) with
$$
y_{j} = \cases 
0, &\text{if $-\ell\le j\le -1$},\\
1, & \text{if $j=0$},\\
(-\beta\gamma^{-1})^{j}\frac{\prod_{i=0}^{j-1}c'_{i}}{\prod_{i=1}^{j}a'_{i}}, &\text{if $1\le j\le \ell$},
\endcases\tag 3.25
$$
where
$$
\frac{\prod_{i=0}^{j-1}c'_{i}}{\prod_{i=1}^{j}a'_{i}} = (-1)^{(\ell+1)j}\bold{i}^{[\frac{\ell+j}{2}]-[\frac{\ell}{2}]+j}q^{j}
\frac{((-q)^{-\ell};-q)^{\frac{1}{2}}_{j}((-q)^{\ell+1};-q)^{\frac{1}{2}}_{j}}{(t-t^{-1})^{j}(q;q)_{j}};
\tag 3.26
$$
or
$$
y_{j} = \cases 
(-\alpha\gamma^{-1})^{-j}\frac{\prod_{i=0}^{-j-1}b'_{-i}}{\prod_{i=1}^{-j}a'_{-i}}, &\text{if $-\ell\le j\le -1$},\\
1, & \text{if $j=0$},\\
0, &\text{if $1\le j\le \ell$}, 
\endcases\tag 3.27
$$
where
$$
\frac{\prod_{i=0}^{-j-1}b'_{-i}}{\prod_{i=1}^{-j}a'_{-i}} = (-1)^{j\ell}\bold{i}^{[\frac{\ell+j}{2}]-[\frac{\ell}{2}]+j^{2}}q^{j}
\frac{((-q)^{-\ell};-q)^{\frac{1}{2}}_{-j}((-q)^{\ell+1};-q)^{\frac{1}{2}}_{-j}}{(q-q^{-1})^{-j}(q;q)_{-j}}.
\tag 3.28
$$
respectively.
\endproclaim
\par
The corresponding left invariant subalgebras are described by the following theorem.
\proclaim{Theorem 3.7} Let $\alpha, \beta$, and $\gamma$ be nonzero complex numbers.
\par
(1) The left invariant subalgebra of $\Cal{A}$ corresponding to $x_{\alpha\beta}$ is spanned by 
$$
\sum_{j\in I_{\ell}}y_{j}m^{(\ell)}_{ij}\sigma^{s},\quad \ell\in \Bbb{Z}_{+},\quad i\in I_{\ell},\quad s=0,1,\tag 3.29
$$
where $y_{j}$ are defined by (3.23) and (3.24).
\par
(2) The left invariant subalgebra of $\Cal{A}$ corresponding to $x_{\alpha\gamma}$ (resp. $x_{\beta\gamma}$) is spanned by
$$
\sum_{j\in I_{\ell}}y_{j}m^{(\ell)}_{ij}\sigma^{s},\quad i\in I_{\ell},\tag 3.30
$$
where $\ell\in \Bbb{Z}_{+}$ and $s=0,1$ satisfy $\ell+s \in 2\Bbb{Z}_{+}$, and $y_{j}$ are defined by (3.25) and (3.26) (resp. (3.27) and (3.28)).
\endproclaim
\par
\head
4. Zonal spherical functions and Askey-Wilson polynomials
\endhead
In this section, we consider the roles of some polynomials from the Askey-Wilson scheme play in the invariants we discussed in the previous section.  We keep the notation of the previous section.  Recall the definitions (see [Ko]) of the $q$-shifted factorials and the $q$-hypergeometric series
$$
\gathered
(a;q)_{m} = \prod_{k=0}^{m-1}(1-aq^{k}), \qquad
(a_{1}, ..., a_{r};q)_{m} = \prod_{j=1}^{r}(a_{j};q)_{m},\\
{}_{s+1}\phi_{s}\left[\matrix a_{1}, ..., a_{s+1}\\ b_{1}, ..., b_{s}\endmatrix ; q, z \right] = \sum_{k=0}^{\infty}\frac{(a_{1}, ..., a_{s+1};q)_{k}z^{k}}{(b_{1}, ...,b_{s}; q)_{k}(q;q)_{k}}.
\endgathered
$$
Recall also the Askey-Wilson polynomials
$$
p_{n}(cos\theta; a, b, c, d|q)=a^{-n}(ab,ac,ad;q)
{}_{4}\phi_{3}\left[\matrix q^{-n},q^{n-1}abcd, ae^{\bold{i}\theta}, ae^{-\bold{i}\theta}\\ ab, ac, ad\endmatrix ; q, q \right],
$$
which satisfy the recurrence relation 
$$
\multline
-(1-q^{-n})(1-q^{n-1}abcd)P_{n}(e^{\bold{i}\theta})\\
=A(\theta) )(P_{n}(qe^{\bold{i}\theta}) - P_{n}(e^{\bold{i}\theta})) + A(-\theta)( P_{n}(q^{-1}e^{\bold{i}\theta}) - P_{n}(e^{\bold{i}\theta})),
\endmultline
$$
and the dual $q$-Krawtchouk polynomials
$$
R_{n}(q^{-x}-q^{x-N-c};q^{c},N|q)= {}_{3}\phi_{2}(q^{-n},q^{-x}, -q^{x-N-c};q^{-N};q,q).
$$
\par
Let $\ell\in \Bbb{Z}_{+}$ and let
$$
v =\sum_{i\in I_{\ell}}x_{i}\xi^{(\ell)}_{i}\sigma^{s}.\tag 4.1
$$ 
The condition $R_{x_{\alpha\beta}}(v)=0$ is equivalent to 
$$
\split
&\alpha\bold{i}^{[\frac{\ell+i-1}{2}]-[\frac{\ell+i}{2}]}t^{\ell-i+1}\frac{(1-t^{-2(\ell-i+1)})^{\frac{1}{2}}(1-t^{-2(\ell+i)})^{\frac{1}{2}}}{q-q^{-1}}x_{i-1}\\
&\qquad -\beta\bold{i}^{[\frac{\ell+i+1}{2}]-[\frac{\ell+i}{2}]}t^{\ell-i-2}\frac{(1-t^{-2(\ell+i+1)})^{\frac{1}{2}}(1-t^{-2(\ell-i)})^{\frac{1}{2}}}{1-t^{-2}}x_{i+1} = 0,
\endsplit\tag 4.2
$$
with the convention that $x_{-\ell-1}= x_{\ell+1}=0$.  We can rewrite (4.2) as
$$
\split
&\bold{i}\alpha\frac{t^{\frac{1}{2}}}{t+t^{-1}}(t^{-\ell+i-1}-t^{\ell-i+1})^{\frac{1}{2}}(t^{-\ell-i}-t^{\ell+i})^{\frac{1}{2}}x_{i-1}\\
&\qquad -\beta t^{-\frac{3}{2}}(t^{-\ell-i-1}-t^{\ell+i+1})^{\frac{1}{2}}(t^{-\ell+i}-t^{\ell-i})^{\frac{1}{2}}x_{i+1} = 0.
\endsplit\tag 4.3
$$
\par
If we choose
$$
\alpha = 1+t^{-2}\qquad \beta = \bold{i}t^{2},\tag 4.4
$$
then (4.3) becomes
$$\split
&\bold{i} t^{-\frac{1}{2}}(t^{-\ell+i-1}-t^{\ell-i+1})^{\frac{1}{2}}(t^{-\ell-i}-t^{\ell+i})^{\frac{1}{2}}x_{i-1}\\
&\qquad - \bold{i}t^{\frac{1}{2}}(t^{-\ell-i-1}-t^{\ell+i+1})^{\frac{1}{2}}(t^{-\ell+i}-t^{\ell-i})^{\frac{1}{2}}x_{i+1} = 0.
\endsplit\tag 4.5
$$
\par
As indicated in [Ko], it is helpful to replace the above homogenous system with the problem of finding the possible eigenvectors of the coefficient matrix of (4.5), i.e., consider the following system of equations
$$\split
&\bold{i} t^{-\frac{1}{2}}(t^{-\ell+i-1}-t^{\ell-i+1})^{\frac{1}{2}}(t^{-\ell-i}-t^{\ell+i})^{\frac{1}{2}}x_{i-1}\\
&\qquad - \bold{i}t^{\frac{1}{2}}(t^{-\ell-i-1}-t^{\ell+i+1})^{\frac{1}{2}}(t^{-\ell+i}-t^{\ell-i})^{\frac{1}{2}}x_{i+1} = \lambda x_{i}.
\endsplit\tag 4.6
$$
The solutions of (4.5) are given by the eigenvectors correspond to the zero eigenvalue of (4.6).  Note that (4.6) is a special case of the identity [Ko, (4.9)], so [Ko, Theorem 4.3] implies that the eigenvalues are
$$
\lambda_{i} = \frac{t^{2i}-t^{-2i}}{t-t^{-1}}, \quad i\in I_{\ell},
$$ 
and the corresponding eigenvectors are given by constant multiples of 
$$
\sum_{n=0}^{2\ell}\bold{i}^{-n}t^{\frac{1}{2}n(n-1)}(t^{2};t^{2})^{-\frac{1}{2}}_{n}(t^{4\ell};t^{-2})^{\frac{1}{2}}_{n}R_{n}(i,2\ell|t^{2})\xi^{(\ell)}_{n-\ell}\sigma^{s},
$$
where 
$$
\split
R_{n}(i,2\ell|t^{2})&= {}_{3}\phi_{2}(t^{-2n}, t^{-2i-2\ell}, -t^{2i-2\ell}; 0,t^{-4\ell}; t^{2}, t^{2})\\
&={}_{3}\phi_{2}((-q)^{-n}, (-q)^{-i-\ell}, -(-q)^{i-\ell}; 0,q^{-2\ell}; -q, -q)
\endsplit
$$
is the dual $q$-Krawtchouk polynomial.  Hence the right invariants of $x_{\alpha\beta}$ (see Theorem 3.5) can also be expressed in terms of some special dual $q$-Krawtchouk polynomials.  Similar discussions hold for the left invariants.
\par
The spherical functions are obtained by taking the intersections of the right and the left invariants.  If the right invariants of $x_{\alpha\beta\gamma}$ in $V^{L}_{\ell}\sigma^{s}$ are given by constant multiples of 
$$
\rho^{(\ell)}_{\alpha\beta\gamma} =\sum_{i\in I_{\ell}}x_{i}\xi^{(\ell)}_{i}\sigma^{s},
$$
and the left invariants of $x_{\alpha'\beta'\gamma'}$ in $V^{R}_{\ell}\sigma^{s}$ are given by constant multiples of
$$
\lambda^{(\ell)}_{\alpha\beta\gamma} =\sum_{j\in I_{\ell}}y_{j}\eta^{(\ell)}_{j}\sigma^{s},
$$
then the $(x_{\alpha\beta\gamma}, x_{\alpha'\beta'\gamma'})$-spherical functions in $M_{\ell}\sigma^{s}$ are given by constant multiples of 
$$
\sum_{i,j\in I_{\ell}}x_{i}y_{j}m^{(\ell)}_{ij}\sigma^{s}.
$$
\par
Consider $(x_{\alpha\beta}, x_{\alpha'\beta'})$-spherical functions. By Proposition 3.4 and Proposition 3.6, there exists a one-dimensional subspace of $(x_{\alpha\beta}, x_{\alpha'\beta'})$-spherical functions in every $M_{\ell}\sigma^{s}$, with $x_{i}$ given by (3.11) and $y_{j}$ given by (3.23).  The spherical functions in $M_{1}\sigma^{s}$ are given by constant multiples of 
$$
\rho_{1}=\left(a^{2} + (\alpha')^{-1}\beta' t(1-q)b^{2} + \alpha\beta^{-1} tq(1-q)^{-1}c^{2} + \alpha\beta^{-1}(\alpha')^{-1}\beta' q^{2}d^{2} \right)\sigma^{s}.\tag 4.7
$$
Note that by (2.7),
$$
\split
k^{m}(\rho_{1})&=(-1)^{ms}(t^{2m}+\alpha\beta^{-1}(\alpha')^{-1}\beta' q^{2}t^{-2m})\\
&=(-1)^{m(s+1)}(q^{m}+\alpha\beta^{-1}(\alpha')^{-1}\beta' q^{-m+2}).
\endsplit\tag 4.8
$$
Since we have the following tensor product decomposition
$$
(M_{1}\sigma^{s})^{\otimes n} = \bigotimes^{n}_{i=0}M_{i}\sigma^{si},\tag 4.9
$$
up to constant multiples, the $(x_{\alpha\beta}, x_{\alpha'\beta'})$-spherical functions in $M_{n}\sigma^{sn}$ are given by a polynomial $P_{n}(\rho_{1})$ of degree n in the variable $\rho_{1}$. To obtain information on this polynomial, consider the following element of $\Cal {U}$
$$
\Gamma = ef + \frac{kt^{-1}+k^{-1}t}{(q-q^{-1})(t-t^{-1})} = -fe + \frac{kt+k^{-1}t^{-1}}{(q-q^{-1})(t-t^{-1})}. \tag 4.10
$$ 
This element belongs to the super center of $\Cal{U}$ (i.e. commutes with even elements and super commutes with odd elements).  By (3.3), (3.11), (3.22) and (3.23), we have
$$
\aligned
R_{\Gamma}(\rho^{\ell}_{\alpha\beta})&=(-1)^{s}\frac{t^{2\ell+1}+t^{-2\ell-1}}{(q-q^{-1})(t-t^{-1})}\rho^{\ell}_{\alpha\beta},\\
L_{\Gamma}(\lambda^{\ell}_{\alpha\beta})&=(-1)^{s}\frac{t^{2\ell+1}+t^{-2\ell-1}}{(q-q^{-1})(t-t^{-1})}\lambda^{\ell}_{\alpha\beta}.
\endaligned\tag 4.11
$$
We claim that 
$$
\split
t(t-t^{-1})(q-q^{-1})k^{m}\Gamma\in &\frac{1-vq^{2m}}{1+vq^{2m-1}}k^{m+1} + \frac{vq^{2m-1}-q}{1+vq^{2m-1}}k^{m-1}\\
& + \Cal{U}x_{\alpha'\beta'} + x_{\alpha\beta}\Cal{U}, 
\endsplit\tag 4.12
$$
where $v=(\alpha')^{-1}\beta'\alpha\beta^{-1}$.
For $x$ and $y$ in $\Cal{U}$, we write $x\sim y$ if $y\in x+ \Cal{U}x_{\alpha'\beta'} + x_{\alpha\beta}\Cal{U}$.  To prove the claim, we first note that
$$
\split
k^{m}ef &= q^{m}(e+\frac{\beta'}{\alpha'}fk)k^{m}f - q^{m-1}\frac{\beta'}{\alpha'}fk^{m}(fk + \frac{\alpha}{\beta}e)+ vq^{m-1}fk^{m}e\\
&\sim vq^{2m-1}k^{m}\frac{k-k^{-1}}{q-q^{-1}}-vq^{2m-1}k^{m}ef,
\endsplit
$$
implies
$$
(1+vq^{2m-1})k^{m}ef \sim vq^{2m-1}\frac{k^{m+1}-k^{m-1}}{q-q^{-1}}. \tag 4.13
$$
Then we add 
$$
(1+vq^{2m-1})k^{m}\frac{kt^{-1}+k^{-1}t}{(q-q^{-1})(t-t^{-1})}
$$
to both sides of (4.13) and multiply through by $t(t-t^{-1})(q-q^{-1})(1+vq^{2m-1})^{-1}$ to prove (4.12). 
\par
By (1.6) and (4.10), we have
$$
\split
(k^{m}\Gamma)(P_{n}(\rho_{1})) &= k^{m}(L_{\Gamma}(P_{n}(\rho_{1})))\\
& = (-1)^{s}\frac{t^{2n+1}+t^{-2n-1}}{(q-q^{-1})(t-t^{-1})} k^{m}(P_{n}(\rho_{1})). \endsplit\tag 4.14
$$
Let 
$$
k^{m}(P_{n}(\rho_{1}))= P_{n}(q^{m}).
$$
Then (4.12) and (4.14) imply that
$$
\split
&((-1)^{s+n}(q^{-m}-q^{m+1})+q-1) P_{n}(q^{m})\\
&\qquad = \frac{1-vq^{2m}}{1+vq^{2m-1}}( P_{n}(q^{m+1}) - P_{n}(q^{m}))\\
&\qquad + \frac{vq^{2m-1}-q}{1+vq^{2m-1}}( P_{n}(q^{m-1}) - P_{n}(q^{m})). 
\endsplit\tag 4.15
$$
\par
Consider the special case when $(\alpha')^{-1}\beta'\alpha\beta^{-1}=1$. Note that
$$
\split
\frac{1-q^{2m}}{1+q^{2m-1}} &= \frac{(1-q^{2m})(1+q^{2m})}{(1+q^{2m-1})(1+q^{2m})}\\
& = \frac{(1-q^{m})(1+q^{m})(1-\bold{i}q^{m})(1+\bold{i}q^{m})}{(1+q^{2m-1})(1+q^{2m})},\endsplit
$$
so if we let $q^{m} = te^{\bold{i}\theta}$, then
$$
\frac{1-q^{2m}}{1+q^{2m-1}} = \frac{(1- te^{\bold{i}\theta })(1+ te^{\bold{i}\theta})(1-\bold{i} te^{\bold{i}\theta})(1+\bold{i} te^{\bold{i}\theta })}{(1-e^{2\bold{i}\theta})(1-qe^{2\bold{i}\theta})}. \tag 4.16
$$
Denote this last expression by $f(\theta)$, then 
$$
f(-\theta) = \frac{q^{2m-1}-q}{1+q^{2m-1}}.
$$
By (4.9), we can replace $s$ in (4.15) by $ns'$.  So we can write (4.15) as
$$
\split
&((-1)^{n(s'+1)}(q^{-m}-q^{m+1})+q-1) P_{n}(q^{m}) \\
&\qquad = f(\theta)( P_{n}(q^{m+1}) - P_{n}(q^{m})) + f(-\theta)( P_{n}(q^{m-1}) - P_{n}(q^{m})). \endsplit\tag 4.17 
$$
\par
If we start with $M_{1}\sigma$ (compare with (4.9)), then $s'=1$, and by (4.8), we can write (4.17) as 
$$
\split
&(q^{-m}-q^{m+1}+q-1) P_{n}(te^{\bold{i}\theta})\\
&\qquad = f(\theta)( P_{n}(qte^{\bold{i}\theta}) - P_{n}(te^{\bold{i}\theta})) + f(-\theta)( P_{n}(q^{-1}te^{\bold{i}\theta}) - P_{n}(te^{\bold{i}\theta})). 
\endsplit\tag 4.18 
$$
Compare (4.18) with the recurrence relation for the Askey-Wilson's polynomials, we have
$$
P_{n}(te^{\bold{i}\theta}) = constant\times {}_{4}\phi_{3}\left[\matrix q^{-n},-q^{n+1}, te^{\bold{i}\theta}, te^{-\bold{i}\theta}\\ q, -\bold{i}q, \bold{i}q \endmatrix ; q, q \right].
$$
Therefore, if $(\alpha')^{-1}\beta'\alpha\beta^{-1}=1$, the $(x_{\alpha\beta}, x_{\alpha'\beta'})$-spherical functions in $M_{n}\sigma^{n}$ are constant multiples of the Askey-Wilson polynomial 
$$
p_{n}(\rho_{1}; t, -t, \bold{i}t, -\bold{i}t | q).
$$
\par
Let us now consider the $(x_{\alpha\beta}, x_{\beta'\gamma'})$-spherical functions as another example.  In this case, we have
$$
\split
k^{m}\Gamma\in & \left(-q^{2m-1}\frac{\beta(\gamma')^{2}}{\alpha(\beta')^{2}} +\frac{t^{-1}}{(q-q^{-1})(t-t^{-1})}\right)k^{m+1}\\
& + \left(-q^{2m-2}\frac{\beta(\gamma')^{2}}{\alpha(\beta')^{2}} +\frac{t}{(q-q^{-1})(t-t^{-1})}\right)k^{m-1}\\
& + q^{2m-1}\frac{\beta(\gamma')^{2}}{\alpha(\beta')^{2}}(1+q^{-1})k^{m} +\Cal{U}x_{\beta'\gamma'} + x_{\alpha\beta}\Cal{U}. 
\endsplit\tag 4.19
$$
If we choose 
$$\alpha=q-q^{-1}, \quad \beta=1, \quad \beta'=t-t^{-1}, \quad \gamma'=1,\tag 4.20
$$
then (4.19) becomes
$$
\split
k^{m}\Gamma &\sim \frac{-q^{2m-1}+1+q^{-1}}{(q-q^{-1})(t-t^{-1})^{2}}(k^{m+1}-k^{m})\\
&\qquad -\frac{q^{2m-2}+1+q}{(q-q^{-1})(t-t^{-1})^{2}}(k^{m-1}-k^{m})\\
&\qquad + \frac{t+t^{-1}}{(q-q^{-1})(t-t^{-1})}k^{m}. 
\endsplit\tag 4.21
$$
For the choice of (4.20), the right $x_{\alpha\beta}$-invariants in $V^{L}_{1}\sigma$ are constant multiples of 
$$
\rho^{1}_{\alpha\beta} = \xi^{(1)}_{-1}\sigma + \bold{i}t(1+q)\xi^{1}_{1}\sigma, \tag 4.22
$$
and the left $x_{\beta'\gamma'}$-invariants in $V^{R}_{1}\sigma$ are constant multiples of 
$$
\lambda^{1}_{\beta'\gamma'} = (1-q)^{\frac{1}{2}}\eta^{(1)}_{0}\sigma + q(1+q)\eta^{1}_{1}\sigma. \tag 4.23
$$
Hence the spherical functions in $M_{1}\sigma$ are constant multiples of 
$$
\rho_{1}=((1-q)^{\frac{1}{2}}ab + \bold{i}q(1+q)b^{2} + t(1+q)(1-q)^{\frac{1}{2}}dc + \bold{i}tq(1+q)^{2}d^{2})\sigma.\tag 4.24
$$
By (4.9), the spherical functions in $M_{n}\sigma^{n}$ are constant multiples of a polynomial $P_{n}(\rho_{1})$ of degree $n$ in $\rho_{1}$.  Let $k^{m}(P_{n}(\rho_{1}))= R_{n}(q^{m})$. Then by (4.21), $ R_{n}(q^{m})$ satisfy the following recurrence relation
$$
\split
&(q^{2m-1}-q^{-1}-1)(R_{n}(q^{m+1})-R_{n}(q^{m}))\\
& + (q^{2m-2}+q+1)(R_{n}(q^{m-1})-R_{n}(q^{m}))\\
&\qquad =(q^{n+1} + q^{n} + q - q^{-n-1} - q^{-n} - q^{-1})R_{n}(q^{m}). 
\endsplit\tag 4.25
$$
\par

\enddocument